\documentclass{owrart}
 \usepackage{graphicx}

\begin{document}


\begin{talk}{Lukas Einkemmer}
{Structure preserving numerical methods for the Vlasov equation}

\noindent
In astro- and plasma physics the behavior of a collisionless plasma is modeled by the Vlasov equation
\[
\partial_{t}f(t,x,v)+v\cdot\nabla_x f(t,x,v)+F\cdot\nabla_{v}f(t,x,v)=0.
\]
The force $F$ is given by the Lorentz force law and the electric $E$ and magnetic $B$ fields are self-consistently determined from the particle-density $f$. Thus, in the most general setting we have to solve the Vlasov equation coupled with Maxwell's equations (the so-called Vlasov--Maxwell system). In many applications, however, it is sufficient to consider the electrostatic case. That is, the Vlasov equation is only coupled to a Poisson problem
\[ \Delta \phi(t,x) = \int f(t,x,v)\,\mathrm{d}v - 1. \]
The potential $\phi$ is then used to determine the electric field by using the relation $E = \nabla_x \phi$. This is the model that we will consider in the present report.

Solving the Vlasov--Poisson system numerically is a challenging task as
\begin{itemize}
	\item the problem is posed in an up to six-dimensional phase space;
	\item the system is nonlinear;
	\item the characteristic time scale in many applications is on the order of picoseconds (the plasma frequency) while interesting physical phenomena can happen on much larger timescales.
\end{itemize}
Here we will focus mostly on the last point. The discrepancy in timescales implies that in order to obtain meaningful results we have to perform a large number of time steps. Therefore, it is not possible to use the numerical scheme in the asymptotic regime. Nevertheless, numerical simulations can still prove useful if they are able to capture the plasma dynamics qualitatively. Thus, in the present context the goal is to construct a numerical integrator that can capture as much of the structure of the analytic solution as possible.

In the case of the Vlasov equation it is well known that an infinite number of conserved quantities exist (for example, all $L^p$ norms are conserved by the analytic solution). Certainly, it is unrealistic to expect that all those invariants are conserved by a numerical approximation. Thus, we focus on the physically important invariants: mass, momentum, energy, and positivity. In addition, we will include entropy and the $L^2$ norm as a measure of how much dissipation is introduced by the numerical scheme.

Since an explicit time stepping scheme needs to obey the CFL constraint given by $v \tau < h$ ($\tau$ is the time step size and $h$ is the grid spacing), the splitting approach introduced by Strang \& Knorr \cite{cheng1976} has been almost universally employed. This approach conserves mass, momentum, all $L^p$ norms, and entropy. Furthermore, an extension to the full Vlasov--Maxwell system has been proposed recently (see \cite{einkemmer1401}).

In addition to the good conservation properties, the splitting approach has another decisive advantage; namely, it reduces the problem of solving the up to six-dimensional Vlasov equation to a sequence of one-dimensional advection equations of the form 
	\[ \partial_t u(t,\xi,\eta) = a(\eta) \partial_{\xi} u(t,\xi,\eta). \]
For this problem it is straightforward to obtain the characteristics analytically and thus semi-Lagrangian methods have become the standard approach. Note, however, that since the feet of the characteristics not necessarily coincide with the grid used, an interpolation procedure has to be employed. Cubic spline interpolation is the most commonly used approach.

More recently, the semi-Lagrangian discontinuous Galerkin (sLdG) scheme has been introduced \cite{qiu2011,rossmanith2011,crouseilles2011}. The main advantage of this approach, compared to spline interpolation (or Fourier based methods), is that the resulting numerical scheme only needs data from at most two adjacent cells in order to compute the advection. This is a particularly important feature if the numerical scheme is implemented on a parallel computer system. A convergence analysis of this method has been conduced in \cite{einkemmer2014convergence,einkemmer2014}.

Some interesting properties of the semi-Lagrangian discontinuous Galerkin approach (sLdG) should be noted.
\begin{itemize}
	\item The error in energy includes additional error terms for the methods of order one and two. Thus, it is prudent to at least use a third order approximation in space.
	\item In all our simulations positivity preservation is less of an issue for the sLdG scheme compared to cubic spline interpolation. In any case, due to the local nature of the numerical method positivity limiters can be easily added (at the cost of some additional diffusion; see \cite{rossmanith2011,qiu2011})
	\item Numerical simulations indicate that while cubic spline interpolation violates the second law of thermodynamics (i.e.~it decreases entropy), the semi-Lagrangian discontinuous Galerkin scheme does not suffer from this deficiency. However, a proof of this statement is still missing (see \cite{einkemmer2016}).
	\item Both for the conserved quantities and for the qualitative features of the solution there is a significant advantage in using higher order approximations (see \cite{einkemmer2016}). This, however, is unexpected based on the regularity of the solution.
\end{itemize}
	To illustrate some of these properties the results for the two-stream instability are shown in Figure \ref{fig:ts}. 

\begin{figure}
	\includegraphics[width=0.70\textwidth]{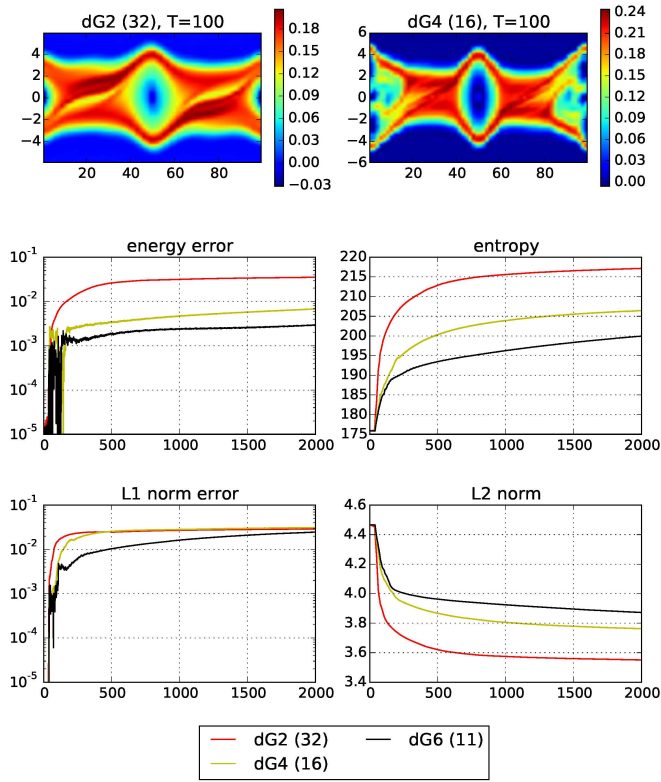}
	\caption{This figure shows $f$ for $t=100$ and the time evolution in the error of the conserved quantities energy, entropy, $L^1$ norm (positivity), and $L^2$ norm. For all numerical schemes $64$
degrees of freedom are employed per space dimension. The order of the discontinuous Galerkin (dG) method is indicated and the number of cells are given in parenthesis.
		\label{fig:ts}}
\end{figure}

\bibliographystyle{plain}
\bibliography{owr-einkemmer.bib}

\begin{thebibliography}{1}

\bibitem{cheng1976}
C.~Cheng and G.~Knorr.
\newblock {The integration of the Vlasov equation in configuration space}.
\newblock {\em J. Comput. Phys.}, 22(3):330--351, 1976.

\bibitem{einkemmer1401}
N.~Crouseilles, L.~Einkemmer, and E.~Faou.
\newblock {A Hamiltonian splitting for the Vlasov--Maxwell system}.
\newblock {\em J. Comput. Phys.}, 238:224--240, 2015.

\bibitem{crouseilles2011}
N.~Crouseilles, M.~Mehrenberger, and F.~Vecil.
\newblock {Discontinuous Galerkin semi-Lagrangian method for Vlasov-Poisson}.
\newblock In {\em ESAIM: Proceedings}, volume~32, pages 211--230. EDP Sciences,
  2011.

\bibitem{einkemmer2016}
L.~Einkemmer.
\newblock {On the geometric properties of the semi-Lagrangian discontinuous
  Galerkin scheme for the Vlasov-Poisson equation}.
\newblock {\em arXiv preprint, arXiv:1601.02280}, 2016.

\bibitem{einkemmer2014}
L.~Einkemmer and A.~Ostermann.
\newblock {Convergence analysis of a discontinuous Galerkin/Strang splitting
  approximation for the Vlasov--Poisson equations}.
\newblock {\em SIAM J. Numer. Anal.}, 52(2):757--778, 2014.

\bibitem{einkemmer2014convergence}
L.~Einkemmer and A.~Ostermann.
\newblock {Convergence analysis of Strang splitting for Vlasov-type equations}.
\newblock {\em SIAM J. Numer. Anal.}, 52(1):140--155, 2014.

\bibitem{qiu2011}
J.M. Qiu and C.W. Shu.
\newblock {Positivity preserving semi-Lagrangian discontinuous Galerkin
  formulation: theoretical analysis and application to the Vlasov--Poisson
  system}.
\newblock {\em J. Comput. Phys.}, 230(23):8386--8409, 2011.

\bibitem{rossmanith2011}
J.A. Rossmanith and D.C. Seal.
\newblock {A positivity-preserving high-order semi-Lagrangian discontinuous
  Galerkin scheme for the Vlasov--Poisson equations}.
\newblock {\em J. Comput. Phys.}, 230(16):6203--6232, 2011.

\end{thebibliography}

\end{talk}

\end{document}